\documentclass[12pt,fleqn]{article}

\usepackage{latexsym}

  \usepackage{amsfonts}
  \usepackage{amsmath}
\usepackage{amssymb}
\usepackage{amscd}
\usepackage{amsthm}

\usepackage{color}
\usepackage[dvips]{graphicx}

\usepackage{authblk}

\newcommand{\dis}{\displaystyle}
\newcommand{\rf}[1]{(\ref{#1})}
\newcommand{\be}{\begin{equation}}
\newcommand{\ee}{\end{equation}}
\newcommand{\ba}{\begin{array}}
\newcommand{\ea}{\end{array}}

\newcommand{\ods}{\par\vspace{0.2cm}\par}

\numberwithin{equation}{section}

\begin{document}

\title{Improving the accuracy of the fast inverse square root algorithm}

\author[1]{\bf Cezary J.\ Walczyk\thanks{walcez@gmail.com}}
\author[2]{\bf  Leonid V. Moroz\thanks{moroz\_lv@polynet.lviv.ua}}
\author[1]{\bf Jan L.\ Cie\'sli\'nski\thanks{j.cieslinski@uwb.edu.pl}} 

\affil[1]{\small Uniwersytet w Bia\l ymstoku, Wydzia{\l} Fizyki, ul.\ Cio\l kowskiego 1L,  15-245 Bia\l ystok, Poland}
\affil[2]{\small   Lviv Polytechnic National University, Department of Security Information and Technology, st. Kn. Romana 1/3, 79000 Lviv, Ukraine }

\date{}

\maketitle

\abstract
We present improved algorithms for fast calculation of  the inverse square root for single-precision floating-point numbers. The algorithms are much more accurate than the famous fast inverse square root algorithm and have the same or similar computational cost. The main idea of our work consists in modifying the Newton-Raphson method and demanding that the maximal error is as small as possible.  The modification can be applied when the distribution of Newton-Raphson corrections is not symmetric (e.g., if they are non-positive functions). 
\ods

\noindent  {\bf Keywords}: floating-point arithmetic; inverse square root; magic constant; Newton-Raphson method

\section{Introduction}

Floating-point arithmetic has became widely used in many applications such as 3D graphics, scientific computing and signal processing  \cite{EL2,cit21,cit13,cit14,cit12}, implemented both in hardware and software \cite{cit22, cit25, cit26,cit23, Vii}. 
Many algorithms can be used to approximate elementary functions, including the inverse square root \cite{EL2,cit21, Vii, EL1, cit18,Hand, cit17, cit19,cit20, Cor, AMC-16,Fog}. All of these algorithms require initial seed to start the approximation. The more accurate is the initial seed, the less iterations is required to compute the function.  In most cases the initial seed is obtained from a look-up table (LUT) which is memory consuming.  The inverse square root function is of particular importance because of it is widely used in 3D computer graphics especially in lightning reflections \cite{Eberly,cit15, cit16}. In this paper we consider an algorithm for computing the inverse square root for single-precision IEEE Standard 754 floating-point numbers (type \textbf{float}) using so the called magic constant instead of LUT \cite{cit1,Lomont,War,Mar}. 
\begin{tabbing}
0000000\=122\=5678\=\kill
\>\textit{1.}\> \textbf{float} InvSqrt(\textbf{float} x)\{ \\
\>\textit{2.}\>\>\textbf{float} halfnumber = 0.5f * x; \\
\>\textit{3.}\>\> \textbf{int} i = *(\textbf{int}*) \&x;\\
\>\textit{4.}\>\> i = R-(i$>>$1);\\
\>\textit{5.}\>\>  x = *(\textbf{float}*)\&i;\\
\>\textit{6.}\>\> x = x*(1.5f-halfnumber*x*x); \\
\>\textit{7.}\>\> x = x*(1.5f-halfnumber*x*x); \\
\>\textit{8.}\>\> \textbf{return} x ;\\
\>\textit{9.}\>\}
\end{tabbing}
The code \textit{InvSqrt}  realizes  a fast algorithm for calculation of the inverse square root. It consists of two main parts. Lines \textit{4} and \textit{5}  produce in a very cheap way a quite good zeroth approximation of the inverse square root of  given positive floating-point number $x$. Lines \textit{6} and \textit{7} apply the Newton-Raphson corrections.

The algorithm  \textit{InvSqrt} has numerous applications, see   \cite{cit4,Janh, cit9, cit7, Schat,cit3, Hann, Li,Hsu, Hsieh,Lv,San,Lin}.  The most important among them is 3D computer graphics, where  normalization of vectors is ubiquitous.  \textit{InvSqrt}  is characterized by a high speed, more that 3 times higher than in computing the inverse square root using library functions.  This property is discussed in detail in \cite{cit10}.  The errors of the fast inverse square root algorithm depend on the choice of the ``magic constant'' $R$.  In several theoretical papers  \cite{Lomont,cit6, cit10,cit11, cit5} (see also the Eberly's monograph \cite{Eberly})  attempts were made to determine analytically the optimal value   of the magic constant (i.e., to minimize errors). In general, this optimal value can depend on the number of iterations, which is a general phenomenon \cite{KM}. The   derivation and comprehensive mathematical description  of all steps of the fast inverse square root algorithm is given in our recent paper \cite{MWHHC}.  It is worthwhile to notice that a magic constant appearing in the Ziv's rounding test \cite{Ziv} has been recently derived in an analytic way as well \cite{DLMT}.

In this paper we develop our analytical approach \cite{MWHHC}  to construct  improved algorithms for computing the inverse square root. The first part of the code will be the same (but we treat $R$ as a free parameter to be determined by minimization of the relative error). The second part of the code will be changed significantly. In fact, we propose a modification of the Newton-Raphson formulae which has the same or similar computational cost but improves the accuracy by several times.  

The magic constant in the \textit{InvSqrt} code serves as a low-cost way of generating a reasonably accurate first approximation of the inverse square root. We are going to insert  other ``magic constants'' into this code in order to increase the accuracy of the algorithm without extra costs.

\section{Algorithm \textit{InvSqrt}. An analytical approach}
\label{theor}

In this section we shortly present main results of \cite{MWHHC}. 
We confine ourselves to positive floating-point numbers
\be \label{float}
   x = (1+m_x) 2^{e_x}
\ee
where $m_x \in [0,1)$ and $e_x$ is an integer. In the case of the IEEE-754 standard, a floating-point number is encoded by
32 bits. The first bit corresponds to a sign (in our case this bit is simply equal  to zero), the next 8 bits correspond to
an exponent $e_x$ and the last 23 bits encodes a mantissa $m_x$. The integer encoded by these 32 bits, denoted by $I_x$, is given by
\be   \label{integ}
    I_x = N_m (B + e_x + m_x) 
\ee
where $N_m = 2^{23}$ and $B=127$ (thus $B+e_x = 1,2,\ldots,254$). The lines \textit{3} and \textit{5} of the \textit{InvSqrt} code interprete a number as an integer \rf{integ} or float \rf{float}, respectively. 
The lines \textit{4}, \textit{6} and \textit{7} of the  code can be written as
\be \ba{l}  \label{3-code}
    I_{y_0} = R - \lfloor I_{x}/2 \rfloor , \\[1ex]
    y_1 = \frac{1}{2} y_0 (3 - y_0^2 x) , \\[1ex]
    y_2 = \frac{1}{2} y_1 (3 - y_1^2 x) .
\ea \ee
The first equation produces, in a surprisingly simple way, a good zeroth approximation $y_0$ of the inverse square root $y = 1/\sqrt{x}$. The next equations can be easily recognized as Newton-Raphson corrections. We point out that the  code \textit{InvSqrt} is invariant with respect to the scaling
\be  \label{invar}
    x \rightarrow \tilde x = 2^{- 2n} x , \qquad  y_k \rightarrow \tilde y_k = 2^n y_k \quad (k=0,1,2) , 
\ee
like the equality $y = 1/\sqrt{x}$ itself.  Therefore, without loss of the generality, we can confine our analysis to the interval
\be
    \tilde A := [1,4) .
\ee
The tilde will denote quantities defined on this interval.  
In \cite{MWHHC} we have shown that the function ${\tilde y}_0$ defined by the first equation of \rf{3-code} can be approximated with a very good accuracy by the piece-wise linear function $\tilde{y}_{00}$ given by
\begin{equation}  \label{y0}
\tilde{y}_{00} (\tilde x, t)  = \left\{   \begin{array}{ll}   \displaystyle 
-\frac{1}{4} \tilde{x} + \frac{3}{4} + \frac{1}{8} t     &     \text{for}  \ \ \tilde{x} \in [1,2 )    \\[2ex]   \displaystyle 
-\frac{1}{8} \tilde{x} + \frac{1}{2} + \frac{1}{8} t    &  \text{for}  \ \   \tilde{x} \in   [2, t)    \\[2ex]   \displaystyle 
-\frac{1}{16} \tilde{x} + \frac{1}{2} + \frac{1}{16} t      &  \text{for} \ \  \tilde{x} \in  [t,4) 
\end{array}  \right. 
\end{equation}
where 
\be
t =  2 + 4 m_R + 2 N_m^{-1} ,
\ee 
and $m_R :=   N_m^{-1} R  -  \lfloor N_m^{-1} R \rfloor$ ($m_R$ is a mantissa of the floating-point numbers corresponding to $R$).  We assumed  $m_R < \frac{1}{2}$ and, as a results, we obtained $e_R = 63$ and $t > 2$.

The only difference between $y_0$ produced by the code \textit{InvSqrt} and $y_{00}$ given by \rf{y0} is the definition of $t$, because $t$ related to the code depends (although in a negligible way) on $x$. Namely,
\be
    |{\tilde y}_{00} - {\tilde y}_{0}|  \leqslant \frac{1}{4} N_m^{-1} = 2^{-25} \approx 2.98 \cdot 10^{-8} .
\ee
Taking into account the invariance \rf{invar}, we obtain
\be
     \left| \frac{y_{00} - y_0}{y_0} \right|  \leqslant 2^{-24} \approx 5.96 \cdot 10^{-8} .
\ee
These estimates do not depend on $t$ (in other words, they do not depend on $R$). 
The relative error of the zeroth approximation \rf{y0} is given by
\be  \label{rel0}
\tilde \delta_0 (\tilde x, t) = \sqrt{\tilde x} \ {\tilde y}_{00} (\tilde x, t) - 1
\ee
This is a continuous function with local maxima at
\begin{equation}   \label{x-max}
\tilde{x}_0^{I}=(6+t)/6,\quad \tilde{x}_0^{II}=(4+t)/3,\quad \tilde{x}_0^{III}=(8+t)/3 , 
\end{equation}
given respectively by
\be \ba{l}
\dis \tilde\delta_0 (\tilde x_0^I,t) = - 1 + \frac{1}{2} \left( 1 + \frac{t}{6} \right)^{3/2} \ , \\[2ex]
\dis \tilde\delta_0 (\tilde x_0^{II},t) = - 1 + 2 \left( \frac{1}{3} \left( 1 + \frac{t}{4} \right) \right)^{3/2} \ , \\[2ex]
\dis \tilde\delta_0 (\tilde x_0^{III},t) = - 1 +  \left( \frac{2}{3} \left( 1 + \frac{t}{8} \right) \right)^{3/2} \ .  
\ea \ee
In order to study global extrema of $\tilde \delta_0 (\tilde x, t)$ we need also boundary values: 
\be \ba{l}
\dis \tilde \delta_0 (1,t)= \dis \tilde \delta_0 (4,t) = \frac{1}{8} \left( t - 4 \right) , \\[2ex]
\dis \tilde \delta_0 (2,t)= \frac{\sqrt{2}}{4}  \left( 1 + \frac{t}{2} \right) - 1 , \\[2ex]
\dis  \tilde \delta_0 (t,t)= \frac{\sqrt{t}}{2} - 1 , \\[2ex]
\ea \ee
which are, in fact, local minima. Taking into account
\be \ba{l}
\dis \tilde \delta_0 (1,t) - \tilde \delta_0 (t,t) = \frac{1}{8} \left( \sqrt{t} - 2 \right)^2 \geqslant 0 \ , \\[2ex]
\dis \tilde \delta_0 (2,t) - \tilde \delta_0 (t,t) = \frac{\sqrt{2}}{8} \left( \sqrt{t} - \sqrt{2} \right)^2 \geqslant 0 , 
\ea \ee
we conclude that
\be  \label{min0}
    \min_{\tilde x \in \tilde A} \tilde \delta_0 (\tilde x, t) = \tilde \delta_0 (t,t) < 0 .
\ee
Because $\tilde\delta_0 (\tilde x_0^{III},t) < 0$ for $t \in (2,4)$, the global maximum is one of the remaining local maxima:
\be  \label{max0}
\max_{\tilde x \in \tilde A} \tilde \delta_0 (\tilde x, t) = \max \{ \tilde\delta_0 (\tilde x_0^I,t), \tilde\delta_0 (\tilde x_0^{II},t) \} .
\ee
Therefore, 
\be
  \max_{x\in \tilde A} | \tilde \delta_0 (\tilde x, t) | = \max \{ | \tilde \delta_0 (t,t) |,  \tilde\delta_0 (\tilde x_0^I,t), \tilde\delta_0 (\tilde x_0^{II},t)  \} .
\ee
In order to minimize this value with respect to $t$, i.e., to find $t_0^{r}$ such that 
\be
  \max_{x\in \tilde A} | \tilde \delta_0 (\tilde x, t_0^{r}) | < \max_{x\in \tilde A} | \tilde \delta_0 (\tilde x, t) | \qquad \text{for} t \neq t_0^{r} , 
\ee
we observe that $| \tilde \delta_0 (t,t) |$ is a decreasing function of $t$, while both maxima ($\tilde\delta_0 (\tilde x_0^I,t)$ and $\tilde\delta_0 (\tilde x_0^{II},t)$) are increasing functions. Therefore, it is sufficient to find $t=t_0^I$ and $t=t_0^{II}$ such that 
\be
    | \tilde \delta_0 (t_0^I, t_0^I) | = \tilde\delta_0 (\tilde x_0^I, t_0^I) \ , \qquad 
| \tilde \delta_0 (t_0^{II}, t_0^{II}) | = \tilde\delta_0 (\tilde x_0^{II},t_0^{II}) ,
\ee
and to choose the greater of these two values. In \cite{MWHHC} we have shown that 
\be  \label{twomax}
    | \tilde \delta_0 (t_0^I, t_0^I) | < | \tilde \delta_0 (t_0^{II},t_0^{II}) | .
\ee
Therefore $t_0^{r} = t_0^{II}$ and 
\be  \label{delta0_max}
  \tilde\delta_{0 \max} := \min_{t \in (2,4)} \left( \max_{x\in \tilde A} | \tilde \delta_0 (\tilde x, t) |   \right)  = | \tilde \delta_0 (t_0^{r},t_0^{r}) | . 
\ee 
The following numerical values result from these calculations \cite{MWHHC}:
\be
  t_0^{r} \approx  3.7309796 , \qquad R_0 = 0x5F37642F , 
\qquad \tilde\delta_{0 \max}  \approx 0.03421281 .
\ee
Newton-Raphson corrections for the zeroth approximation given by $\tilde{y}_{00}$ will be denoted by $\tilde{y}_{0k}$ ($k=1,2,\ldots$), in particular: 
\be \ba{l}  \label{NR}
\tilde{y}_{01}(\tilde{x},t) = \frac{1}{2}  \tilde{y}_{00}(\tilde{x},t)(3-\tilde{y}_{00}^2(\tilde{x},t) \,\tilde{x}),   \\[2ex]
\tilde{y}_{02}(\tilde{x},t) =\frac{1}{2} \tilde{y}_{01}(\tilde{x},t)(3-\tilde{y}_{01}^2(\tilde{x},t)\, \tilde{x}). 
\ea \ee
and the corresponding relative error functions will be denoted by $\tilde{\delta}_{k}(\tilde{x},t)$: 
\begin{equation}     \label{delta k}
\tilde{\delta}_{k}(\tilde{x},t):= \frac{\tilde{y}_{0k}(\tilde{x},t) -\tilde y}{\tilde y} =\sqrt{\tilde{x}}\tilde{y}_{0k}(\tilde{x},t)-1, \qquad  (k=0,1,2,\ldots) ,   
\end{equation}
where we included also the case $k=0$, see \rf{rel0}. The obtained approximations of the inverse square root depend on the parameter $t$  directly related to the magic constant $R$.  The value of this parameter can be estimated  by analysing the relative error of  $\tilde{y}_{0k}(\tilde{x},t)$ with respect to $1/\sqrt{\tilde{x}}$. 
As the best estimation we consider  $t=t_k^{(r)}$  minimizing  the relative error  $\tilde{\delta}_{k}(\tilde{x},t)$:
\begin{equation}   \label{r4_2}
\forall_{t\neq t_k^{(r)}}  \left( \tilde\delta_{k \max} \equiv \max_{\tilde{x}\in\tilde{A}}|\tilde{\delta}_{k}(\tilde{x},t_k^{(r)})|<\max_{\tilde{x}\in\tilde{A}}| \tilde{\delta}_{k}(\tilde{x},t)|  \right)   .
\end{equation}
We point out that in general the optimum value of the magic constant can depend on the number of Newton-Raphson corrections. Calculations carried out in \cite{MWHHC} gave the following results: 
\be  \ba{l}  \label{rezultaty}
t_1^{r} = t_2^{r} =3.7298003 , \qquad R_1^{r} = R_2^{r} = 0x5F375A86 , \\[1ex]  
\tilde\delta_{1 \max}  \approx 1.75118 \cdot 10^{-3} , \quad \tilde\delta_{2 \max}  \approx 4.60 \cdot 10^{-6} . 
\ea \ee
We omit details of the computations except one important point. Using \rf{delta k} for expressing  $\tilde{y}_{0k}$ by $\tilde{\delta}_k$ and $\sqrt{\tilde{x}}$ we can  rewrite \rf{NR} as 
\begin{equation}   \label{r4_15}
\tilde{\delta}_k(\tilde{x},t)=-\frac{1}{2}\tilde{\delta}_{k-1}^2 (\tilde{x},t) (3+\tilde{\delta}_{k-1}(\tilde{x},t)), \quad (k=1,2,\ldots) .
\end{equation}
The quadratic dependence on $\tilde{\delta}_{k-1}$ means that every Newton-Raphson correction improves the accuracy by several orders of magnitude, compare \rf{rezultaty}. 

The formula \rf{r4_15} suggests another way of improving the accuracy because the   functions $\tilde{\delta}_{k}$  are always non-positive for any $k \geqslant 1$. Roughly saying, we are going to shift the graph of $\tilde \delta_k$ upwards by an appropriate modification of the Newton-Raphson formula. In the next sections we describe the general idea of this modification and derive two new codes for fast and accurate computation of the inverse square root.

\section{Modified Newton-Raphson formula}
\label{alg0}

The formula \rf{r4_15} shows that Newton-Raphson corrections are nonpositive (see also Fig.~4 and Fig.~5 in \cite{MWHHC}), i.e., they take values in intervals     $[-{\tilde{\delta}}_{k \max},\,0 ]$, where $k=1,2,\ldots$. Therefore, it is natural to introduce a correction term into Newton-Raphson formulas \rf{NR}. We expect that the corrections will be roughly half of the maximal relative error.  Instead of the maximal error we introduce two parameters, $d_1$ and $d_2$. Thus we get modified Newton-Raphson formulas: 
\be \ba{l}   \label{mNR} 
\dis \tilde{y}_{11}(\tilde{x},t,d_1) =2^{-1} \tilde{y}_{00}(\tilde{x},t)(3-\tilde{y}_{00}^2(\tilde{x},t) \,\tilde{x})+\frac{d_1}{2\sqrt{\tilde{x}}}, \\[2ex]
\dis \tilde{y}_{12}(\tilde{x},t,d_1,d_2) =2^{-1} \tilde{y}_{11}(\tilde{x},t,d_1)(3-\tilde{y}_{11}^2(\tilde{x},t,d_1)\, \tilde{x})+\frac{d_2}{2\sqrt{\tilde{x}}},
\ea \ee
where we still assume the zeroth approximation in the form \rf{y0}. 
The corresponding error functions, 
\be
\tilde{\delta}^{''}_{k}(\tilde{x},t,d_1,\ldots,d_k)=\sqrt{\tilde{x}}\, \tilde{y}_{1k}(\tilde{x},t,d_1,\ldots,d_k)-1 , \qquad  k\in\{0,1,2,\ldots \} ,
\ee
(where $\tilde{y}_{10}(\tilde{x},t) := \tilde{y}_{00}(\tilde{x},t)$), satisfy
\begin{equation}  \label{r5_2}
\tilde{\delta}^{''}_k =-\frac{1}{2}\tilde{\delta}^{''2}_{k-1}(3+\tilde{\delta}^{''}_{k-1}) + \frac{d_k}{2},
\end{equation}
where: $\tilde{\delta}^{''}_0(\tilde{x},t)=\tilde{\delta}_0(\tilde{x},t)$. Note that
\be  \label{delty}
\tilde{\delta}^{''}_1(\tilde{x},t,d_1) = \tilde{\delta}_1 (\tilde{x},t) + \frac{1}{2} d_1  . 
\ee
In order to simplify notation we usually will supress the explicit dependence on $d_j$. We will write, for instance, $\tilde{\delta}^{''}_{2}(\tilde{x},t)$ instead of  $\tilde{\delta}^{''}_{2}(\tilde{x},t,d_1,d_2)$. 

The corrections of the form \rf{mNR} will decrease relative errors in comparison with the results of  earlier  papers \cite{Lomont,MWHHC}. We have 3 free parameters ($d_1, d_2$ and $t$) to be determined  by minimizing the maximal error (in principle the new parameterization can give a new estimation of the parameter $t$).  By analogy to (\ref{r4_2}), we are going to find $t=t^{(0)}$  minimizing the error of the first correction (\ref{r4_2}): 
\begin{equation}
\forall_{t\neq t^{(0)}}\max_{\tilde{x}\in\tilde{A}}|\tilde{\delta}_{1}^{''}(\tilde{x},t^{(0)})|<\max_{\tilde{x}\in\tilde{A}}| \tilde{\delta}^{''}_{1}(\tilde{x},t)| ,    \label{r5_3a}
\end{equation}
where, as usual, $\tilde A = [1,4]$.

The first of Eqs.~\rf{r5_2} implies that for any $t$ the maximal value of $\tilde{\delta}^{''}_1 (\tilde x, t)$ equals $\frac{1}{2} d_1$ and is attained at zeros of $\tilde{\delta}^{''}_0 (\tilde x, t)$. Using results of section~\ref{theor}, including \rf{min0}, \rf{max0}, \rf{twomax} and \rf{delta0_max}, we  conclude that the minimum value of $\tilde{\delta}^{''}_1 (\tilde x, t)$ is attained either for $\tilde x = t$ or for $\tilde x = x_0^{II}$ (where there is the second maximum of $\tilde{\delta}^{''}_0 (\tilde x, t)$), i.e.,
\be
  \min_{\tilde x \in \tilde A} \tilde{\delta}^{''}_1 (\tilde x, t)  = \min \left\{   \tilde{\delta}^{''}_1 (t, t), \tilde{\delta}^{''}_1 (x_0^{II}, t)  \right\}
\ee
Minimization of $|\tilde{\delta}^{''}_1(\tilde{x},t)|$ can be done with respect to $t$ and with respect to $d_1$ (these both operations obviously commute). corresponds to 
\begin{equation}   \label{r5_3}
\underbrace{\max_{\tilde{x}\in\tilde{A}} \tilde{\delta}^{''}_{1}(\tilde{x},t^{(0)})}_{\tilde{\delta}^{''}_{1 \max}}=-\min_{\tilde{x}\in\tilde{A}} \tilde{\delta}^{''}_{1}(\tilde{x},t^{(0)}) . 
\end{equation}
Taking into account 
\begin{equation}
\max_{\tilde{x}\in\tilde{A}} \tilde{\delta}^{''}_{1}(\tilde{x},t^{(0)})=\frac{d_1}{2},\quad \min_{\tilde{x}\in\tilde{A}} \tilde{\delta}^{''}_{1}(\tilde{x},t^{(0)})=\tilde{\delta}^{''}_{1}(t^{(0)},t^{(0)})=-\tilde{\delta}_{1\max}+\frac{d_1}{2},\label{r5_3_2}
\end{equation}
we get from \rf{r5_3}: 
\begin{equation}  \label{8,7559}
\tilde{\delta}^{''}_{1 \max} = \frac{1}{2} d_1 =\frac{1}{2}\tilde{\delta}_{1\, \max}\simeq 8.7559\cdot 10^{-4},
\end{equation}
where 
\be  \label{delta1_max}
  \tilde\delta_{1 \max} := \min_{t \in (2,4)} \left( \max_{x\in \tilde A} | \tilde \delta_1 (\tilde x, t) |   \right)   . 
\ee 
and the numerical value of $\tilde{\delta}_{1 \max}$ is given by \rf{rezultaty}.
These conditions are satisfied for
\begin{equation}
t^{(0)}=t_1^{(r)}\simeq 3.7298003 .
\end{equation}

In order to minimize the relative error of the second correction we use equation analogous to \rf{r5_3}: 
\begin{equation}
\underbrace{\max_{\tilde{x}\in\tilde{A}} \tilde{\delta}^{''}_{2}(\tilde{x},t^{(0)})}_{\tilde{\delta}^{''}_{2\max}}=-\min_{\tilde{x}\in\tilde{A}} \tilde{\delta}^{''}_{2}(\tilde{x},t^{(0)}),\label{r5_3_4}
\end{equation}
where from \rf{r5_2} we have
\begin{equation}
\max_{\tilde{x}\in\tilde{A}} \tilde{\delta}^{''}_{2}(\tilde{x},t^{(0)})=\frac{d_2}{2},\quad \min_{\tilde{x}\in\tilde{A}} \tilde{\delta}^{''}_{2}(\tilde{x},t^{(0)})=-\frac{1}{2}\tilde{\delta}^{''2}_{1\max} \left(3+\tilde{\delta}^{''}_{1\max}\right)+\frac{d_2}{2}.\label{r5_3_5}
\end{equation}
Hence 
\begin{equation}
\tilde{\delta}^{''}_{2\,\max}=\frac{1}{4}\tilde{\delta}^{''2}_{1\max} \left(3+\tilde{\delta}^{''}_{1\max}\right) . 
\end{equation}
Expressing this result in terms of formerly computed $\tilde{\delta}_{1\max}$ and $\tilde{\delta}_{2\max}$, we obtain 
\begin{equation}  \label{8 razy}
\tilde{\delta}^{''}_{2\,\max}=\frac{1}{8}\tilde{\delta}_{2\, \max}+\frac{3}{32}\tilde{\delta}_{1\, \max}^3  
\simeq 5.75164\cdot 10^{-7}  \simeq \frac{\tilde{\delta}_{2\, \max}}{7.99}  ,  \end{equation}
where 
\[
\tilde{\delta}_{2\max} =\frac{1}{2}\tilde{\delta}^{2}_{1\max}(3-\tilde{\delta}_{1\max}) .
\]
Therefore, the modification of Newton-Raphson formulas decreased the relative error almost 8 times.

In principle, the presented idea can be applied in any case in order to improved the accuracy of the Newton-Raphson corrections. However, in order to implement this idea in the form of a computer code, we have to replace the unknown $1/\sqrt{\tilde x}$ (i.e., $\tilde y$ in the general case) on the right-hand sides of \rf{mNR} by some approximation. In the case of the inverse square root function this can be done without difficulties. The most natural choice is to replace the unknown inverse square root by its forward or backward approximation. In sections~\ref{alg1} and~\ref{alg2} we present two algorithms resulting from two simplest approximations of $\tilde y$.

\section{A possibility of further minimization of the second correction}

Considerations of the previous section assumed that we fix $d_1$ by minimizing the first Newton-Raphson correction, and then we obtain the optimum value of $d_2$. 
The error of the second correction has the same value for all $t$ from some neighbourhood of $t^{(0)}$:
\begin{equation}  \label{przedzial}
t\in[t^{(0)}_1,t^{(0)}_2]=[3.72978085,3.72981263] ,
\end{equation}
where the boundaries ($t^{(0)}_1$ and $t^{(0)}_2$) are computed as solutions of the  equations: 
\begin{align}
-\tilde{\delta}^{''}_{2\max}=\tilde{\delta}^{''}_{2}(t^{(0)}_1,t^{(0)}_1),\quad -\tilde{\delta}^{''}_{2\max}=\tilde{\delta}^{''}_{2}((4+t^{(0)}_2)/3,t^{(0)}_2) .
\end{align}
The interval \rf{przedzial} corresponds to the following set of magic constants: 
\begin{align}
\{1597463133,\,1597463134,\,\dots,\, 1597463200\}=&\nonumber\\
=\{\mathrm{0x5F375A5D},\,\mathrm{0x5F375A5E},\,...,\,\mathrm{0x5F375AA0}\}&. \label{r5_3_7}
\end{align}
This nonuniqueness is due to the fact that choosing the parameter 
$d_1$ as minimizing the error of the first correction does not lead to the full minimization of the second correction.  

We relax the assumption about minimizing the error of the first correction.  We are going to find the minimum value of the second correction for $d_1^{'}$ which yields  the following minimum and maximum of the first correction: 
\[
\tilde{\delta}_0(t^{(0)},t^{(0)})+d_1^{'}/2\quad\mathrm{i}\quad d_1^{'}/2 \ . 
\]
Then, the condition for minimization of the second correction error is given by:
\be \ba{l} \dis
-\frac{1}{8}d_1^{'2}\left(3+\frac{d_1^{'}}{2}\right)+\frac{d_2^{'}}{2}= \\[3ex]\dis
-\frac{1}{2}\left(\tilde{\delta}_1(t^{(0)},t^{(0)})+\frac{d_1^{'}}{2}\right)^2 \left(3+\underbrace{\tilde{\delta}_1(t^{(0)},t^{(0)})}_{-\tilde{\delta}_{1\max}}+\frac{d_1^{'}}{2}\right)+\frac{d_2^{'}}{2}.
\ea \ee
This is a quadratic equation for $d_1^{'}$:
\[
d_1^{'2}-2d_1^{'}(\tilde{\delta}_{1\max}-2)-\frac{4}{3}\tilde{\delta}_{1\max}(3-\tilde{\delta}_{1\max})=0,
\]
and its positve root is, indeed, lees than $d_1$: 
\[
d_1^{'}=\tilde{\delta}_{1\max}-2+2\sqrt{1-\tilde{\delta}_{1\max}/12}<d_1.
\]
Thus the maximum value of the first correction relative error,  
\[
d_1^{'}/2\simeq 8.75464\cdot 10^{-4}<\tilde{\delta}_{1\max}^{''} \simeq 8.7559 \cdot 10^{-4} \ ,
\]
where $\tilde{\delta}_{1\max}^{''}$ is given by \rf{8,7559}, is less than the modulus of its minimum:
\[
d_1^{'}/2+2-2\sqrt{1-\tilde{\delta}_{1\max}/12}\simeq 8.75720\cdot 10^{-4}>\tilde{\delta}_{1\max}^{''}.
\]
The minimization of the maximal relative error of the second correction reduces to the solution of the following equation, 
\begin{equation}
\frac{d_2^{'}}{2}=\frac{1}{8}d_1^{'2}\left(3+\frac{1}{2}d_1^{'}\right)-\frac{d_2^{'}}{2} ,
\end{equation}
where the maximal error of the second correction is equated with the modulus of the minimal error. Hence
\begin{equation}
\frac{d_2^{'}}{2}=\frac{1}{16}d_1^{'2}\left(3+\frac{1}{2}d_1^{'}\right)\simeq 5.74996\cdot 10^{-7} \ .
\end{equation}
The obtained result is less than $\tilde{\delta}_{2\max}^{''}$, given by \rf{8 razy}, by only  $1.68\cdot 10^{-10}$ which is negligible compared to round-off errors, of order $6\cdot 10^{-8}$, appearing during calculations with the precision \textbf{float}.

The possibility described in this section will not be used in the construction of our algorithms because of at least two reasons. First, the improvement of the accuracy of the second correction is infinitesimal. Second, our aim is to build an algorithm which can be stopped either after one or after two iterations. The best algorithm  for two iterations (discussed in this section) is not optimal when stopped after the first iteration.

\section{Algorithm {\it InvSqrt1}}
\label{alg1}

\begin{figure}
\begin{center}
\includegraphics[width=15cm]{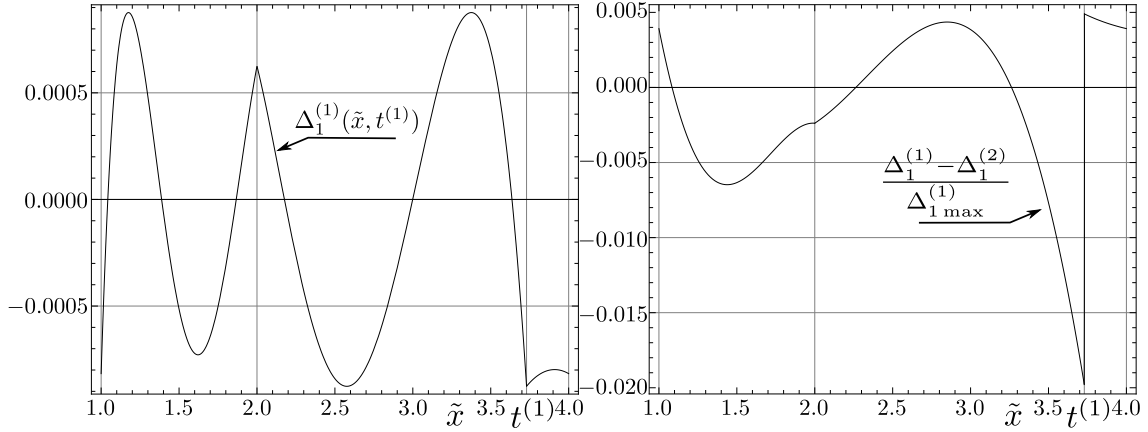}
\end{center}
\caption{Left: graph of the function $\Delta_1^{(1)}(\tilde{x},t^{(1)})$. Right:   $(\Delta_1^{(1)}(\tilde{x},t^{(1)}) - \Delta_1^{(2)}(\tilde{x},t^{(2)}))$ as a fraction of $\Delta_{1\max}^{(1)}$.}
\label{pic6}
\end{figure}

Approximating $1/\sqrt{\tilde{x}}$ by $\tilde{y}_{21}$ and $\tilde{y}_{22}$, respectively, we transform \rf{mNR} into 
\be \ba{l}   \label{mNR1} 
\dis \tilde{y}_{21} = \frac{1}{2} \tilde{y}_{20} (3-\tilde{y}_{20}^2 \,\tilde{x})+ \frac{1}{2} d_1 \tilde{y}_{21} , \\[2ex]
\dis \tilde{y}_{22} = \frac{1}{2} \tilde{y}_{21} (3-\tilde{y}_{21}^2 \, \tilde{x})+\frac{1}{2} d_2 \tilde{y}_{22} ,
\ea \ee
where $\tilde y_{2k}$ ($k=1,2,\ldots$) depend on $\tilde{x},t$ and $d_j$ (for $1 \leqslant j  \leqslant k$). We assume $\tilde y_{20} \equiv \tilde y_{00}$. Thus $\tilde y_{21}$ and $\tilde y_{22}$ can be explicitly expressed by $\tilde y_{20}$ and $\tilde y_{21}$, respectively. The error functions are defined in the usual way: 
\begin{equation}  \label{r5:8}
\Delta_k^{(1)} = \frac{\tilde{y}_{2k} - \tilde{y}}{\tilde{y}} = \sqrt{\tilde{x}}\,\tilde{y}_{2k} - 1 \ . 
\end{equation}
Substituting \rf{r5:8} into \rf{mNR1} we get:
\be  \label{r5:9a}
\Delta_1^{(1)}(\tilde{x},t,d_1)=\frac{d_1}{2-d_1}- \frac{1}{2-d_1}\tilde{\delta}_{0}^2(\tilde{x},t)(3+\tilde{\delta}_{0}(\tilde{x},t)) =\frac{d_1 +2\tilde{\delta}_{1}(\tilde{x},t)}{2-d_1}, 
\ee 
\be  \label{r5:9b}
\Delta_2^{(1)}(\tilde{x},t,d_,d_2)=\frac{d_2}{2-d_2}- \frac{1}{2-d_2}\left( \Delta_{1}^{(1)}(\tilde{x},t,d_1) \right)^2 \left(3+\Delta_{1}^{(1)}(\tilde{x},t,d_1) \right).    
\ee
The equation \rf{r5:9a} expresses $\Delta_1^{(1)}(\tilde{x},t,d_1)$ as a linear function of the nonpositive function $\tilde{\delta}_1(\tilde{x},t)$ with coefficients depending on the parameter $d_1$.  The optimum parameters $t$ and $d_1$ will be estimated by the procedure described in section~\ref{alg0}. First, we minimize the amplitude of the relative error function, i.e., we find $t^{(1)}$ such that 
\begin{equation}   \label{r5_100}
 \max_{\tilde{x}\in\tilde{A}} {\Delta}_{1}^{(1)}(\tilde{x},t^{(1)}) - \min_{\tilde{x}\in\tilde{A}} {\Delta}_{1}^{(1)}(\tilde{x},t^{(1)}) \leqslant \max_{\tilde{x}\in\tilde{A}}  {\Delta}_{1}^{(1)}(\tilde{x},t) - \min_{\tilde{x}\in\tilde{A}}  {\Delta}_{1}^{(1)}(\tilde{x},t) 
\end{equation}
for all $t\neq t^{(1)}$. Second,  we determine $d_1^{(1)}$ such that
\be
    \max_{\tilde{x}\in\tilde{A}} {\Delta}_{1}^{(1)}(\tilde{x},t^{(1)},d_1^{(1)}) = - \min_{\tilde{x}\in\tilde{A}}{\Delta}_{1}^{(1)}(\tilde{x},t^{(1)},d_1^{(1)}) \  .
\ee  
Thus we have
\begin{equation}   \label{r5_10}
\max_{\tilde{x}\in\tilde{A}}|{\Delta}_{1}^{(1)}(\tilde{x},t^{(1)},d_1^{(1)})|\leqslant \max_{\tilde{x}\in\tilde{A}}| {\Delta}_{1}^{(1)}(\tilde{x},t,d_1)| 
\end{equation}
for all real $d_1$ and $t \in (2,4)$.  $\Delta_1^{(1)}(\tilde{x},t)$ is an increasing function of $\tilde{\delta}_1(\tilde{x},t)$, hence  
\begin{equation}
-\frac{d_1^{(1)}-2\max_{\tilde{x}\in\tilde{A}}|\tilde{\delta}_{1}(\tilde{x},t_1^{(1)})|}{2-d_1^{(1)}} =\frac{d_1^{(1)}}{2-d_1^{(1)}},
\end{equation}
which is satisfied for 
\begin{equation}
d_1^{(1)}=\max_{\tilde{x}\in\tilde{A}}|\tilde{\delta}_{1}(\tilde{x},t_1^{(1)})|=\tilde{\delta}_{1\max}.
\end{equation}
Thus we can find the maximum error of the first correction $\Delta_1^{(1)}(\tilde{x},t_1^{(1)})$ (presented in the left part of Fig.~\ref{pic6}):
\begin{equation}
\max_{\tilde{x}\in\tilde{A}}|\Delta_1^{(1)}(\tilde{x},t^{(1)})|= \frac{\max_{\tilde{x}\in\tilde{A}}|\tilde{\delta}_{1}(\tilde{x},t^{(1)})|}{2- \max_{\tilde{x}\in\tilde{A}}|\tilde{\delta}_{1}(\tilde{x},t^{(1)})|},
\end{equation}
which assumes the minimum value for $t^{(1)}=t_1^{(r)}$:
\begin{equation}
\Delta_{1\max}^{(1)}=\frac{\max_{\tilde{x}\in\tilde{A}}|\tilde{\delta}_{1}(\tilde{x},t_1^{(r)})|}{2- \max_{\tilde{x}\in\tilde{A}}|\tilde{\delta}_{1}(\tilde{x},t_1^{(r)})|}= \frac{\tilde{\delta}_{1\,\max}}{2-\tilde{\delta}_{1\,\max}}     \simeq 8.7636\cdot 10^{-4}  \simeq \frac{\tilde{\delta}_{1\,\max}}{2.00} .  \ \ 
\end{equation}

\noindent  Analogically we can determine the value of $d_2^{(1)}$ (assuming that $t=t^{(1)}$ is  fixed): 
\begin{equation}
-\frac{d_2^{(1)}-\max_{\tilde{x}\in\tilde{A}}|\Delta_{1}^{(1)2}(\tilde{x},t^{(1)})(3+\Delta_{1}^{(1)}(\tilde{x},t^{(1)}))|}{2-d_2^{(1)}} =\frac{d_2^{(1)}}{2-d_2^{(1)}}.
\end{equation}
Now, the deepest minimum comes from the global maximum  
\begin{equation}
\max_{\tilde{x}\in\tilde{A}}|\Delta_{1}^{(1)2}(\tilde{x},t^{(1)})(3+\Delta_{1}^{(1)}(\tilde{x},t^{(1)}))|= \frac{2\tilde{\delta}^2_{1\,\max}(3-\tilde{\delta}_{1\,\max})}{(2-\tilde{\delta}_{1\,\max})^3} .
\end{equation}
Therefore we get
\begin{equation}
d_2^{(1)}= \frac{\tilde{\delta}^2_{1\,\max}(3-\tilde{\delta}_{1\,\max})}{(2-\tilde{\delta}_{1\,\max})^3}\simeq 1.15234\cdot 10^{-6},
\end{equation}
and the maximum error of the second correction is given by 
\begin{equation}
\Delta_{2\max}^{(1)}=\frac{d^{(1)}_{2}}{2-d_2^{(1)}}   \simeq 5.76173\cdot 10^{-7} \simeq \frac{\tilde{\delta}_{2\,\max}}{7.98} .
\end{equation}
The computer code for the algorithm described above is the following modification of  $InvSqrt$:
\begin{tabbing}
0000000\=122\=5678\=\kill
\>\textit{1.}\> \textbf{float} InvSqrt1(\textbf{float} x)\{ \\
\>\textit{2.}\>\>\textbf{float} simhalfnumber = 0.50043818f * x; \\
\>\textit{3.}\>\> \textbf{int} i = *(\textbf{int}*) \&x;\\
\>\textit{4.}\>\> i = 0x5F375A86 - (i$>>$1);\\
\>\textit{5.}\>\>  x = *(\textbf{float}*)\&i;\\
\>\textit{6.}\>\> x = x*(1.5013145f - simhalfnumber*x*x); \\
\>\textit{7.}\>\> x = x*(1.5000008f - 0.99912498f*simhalfnumber*x*x); \\
\>\textit{8.}\>\> \textbf{return} x ;\\
\>\textit{9.}\>\}
\end{tabbing}
where
\begin{align*}
0.50043818f &\simeq(2-d_1^{(1)})^{-1},\quad&\quad 1.5013145f &\simeq 3(2-d_1^{(1)})^{-1},\\
1.5000008f &\simeq 3(2-d_2^{(1)})^{-1},\quad&\quad 0.99912498f &\simeq (2-d_1^{(1)})(2-d_2^{(1)})^{-1}.
\end{align*}
Comparing $InvSqrt1$ with $InvSqrt$ we easily see that the number of algebraic operations in $InvSqrt1$ is greater just by 1 (an additional multiplication in line $7$, corresponding to the second iteration of the modified Newton-Raphson procedure).

\section{Algorithm {\it InvSqrt2}}
\label{alg2}

Approximating $1/\sqrt{\tilde{x}}$ by $\tilde{y}_{30}$ and $\tilde{y}_{31}$, respectively, we transfom \rf{mNR} into 
\be \ba{l}   \label{mNR2} 
\dis \tilde{y}_{31} = \frac{1}{2} \tilde{y}_{30} (3-\tilde{y}_{30}^2 \,\tilde{x})+ \frac{1}{2} d_1 \tilde{y}_{30} , \\[2ex]
\dis \tilde{y}_{32} = \frac{1}{2} \tilde{y}_{31} (3-\tilde{y}_{31}^2 \, \tilde{x})+\frac{1}{2} d_2 \tilde{y}_{31} ,
\ea \ee
where $\tilde y_{3k}$ ($k=1,2,\ldots$) depend on $\tilde{x},t$ and $d_j$ (for $j \leqslant k$). We assume $\tilde y_{30} \equiv \tilde y_{00}$. Thus $\tilde y_{31}$ and $\tilde y_{32}$ can be explicitly expressed by $\tilde y_{30}$ and $\tilde y_{31}$, respectively. The error functions are defined in the usual way: 
\begin{equation}  \label{r5:20}
\Delta_k^{(2)} =\sqrt{\tilde{x}}\,\tilde{y}_{3k} - 1 \ . 
\end{equation}
Substituting \rf{r5:20} into \rf{mNR2} we get:
\be  \label{r5:21a}
\Delta_1^{(2)}(\tilde{x},t,d_1)=\frac{1}{2} d_1 \left(1+\tilde{\delta}_0 (\tilde{x},t) \right)- \frac{1}{2}\tilde{\delta}_{0}^2(\tilde{x},t)\left(3+\tilde{\delta}_{0}(\tilde{x},t) \right) 
\ee 
\be  \label{r5:21b}
\Delta_2^{(2)}(\tilde{x},t,d_1,d_2)=\frac{1}{2} d_2 \left(1+{\Delta^{(2)}_1} \right) - \frac{1}{2} \left( \Delta_{1}^{(2)} \right)^2 \left(3+\Delta_{1}^{(2)} \right) ,  
\ee
where $\Delta_{1}^{(2)} \equiv \Delta_{1}^{(2)}(\tilde{x},t,d_1)$.

First, we are going to determine $t$ and $d_1^{(2)}$ minimizing the maximum absolute value of the relative error of the first correction. Therefore, we have to solve the following equation:
\begin{equation}
0=\frac{\partial\Delta_1^{(2)}(\tilde{x},t)}{\partial \tilde{\delta}_{0}(\tilde{x},t)}=\frac{1}{2}d_1^{(2)}-3\tilde{\delta}_{0}(\tilde{x},t) -\frac{3}{2}\tilde{\delta}_{0}^2(\tilde{x},t).    \label{r5:22}
\end{equation}
Its solution
\begin{equation}
\tilde{\delta}^{+}=\sqrt{1+d_1^{(2)}/3}-1\label{r5:23}
\end{equation} 
corresponds to the value 
\begin{equation}
{\Delta}^{(2)}_{1\max}=\frac{1}{2}d_1^{(2)}(1+\tilde{\delta}^{+})- \frac{1}{2}\tilde{\delta}^{+2}(3+\tilde{\delta}^{+}),\label{r5:24}
\end{equation} 
which is a maximum of  $\Delta_1^{(2)}(\tilde{x},t)$ because its second derivative with respect to $\tilde{x}$, i.e.,
\begin{equation}
\partial_{\tilde{x}}^2 \Delta_1^{(2)}(\tilde{x},t)=\partial^2_{\tilde{x}}\tilde{\delta}_0(\tilde{x},t) \partial_{\tilde{\delta}_0}{\Delta}_1^{(2)}(\tilde{x},t)-3(1+\tilde{\delta}_0(\tilde{x},t)) (\partial_x \tilde{\delta}_0(\tilde{x},t))^2,\label{r5:25}
\end{equation}
is negative. In order to determine the dependence of $d_1^{(2)}$ on the parameter  $t$ we solve the equation
\begin{equation}
-\Delta_1^{(2)}(t,t)={\Delta}^{(2)}_{1\max},              \label{r5:26}
\end{equation}
which equates (for some $t=t^{(2)}$) the maximum value of error with the modulus of the minimum value of error. Thus we obtain the following relations:  
\begin{align}
\delta^{+} &=-1-\frac{1}{4}\sqrt{t}+\frac{1}{8}\frac{t}{f(t)}+\frac{1}{2}f(t),\label{r5:27a}\\
d_1^{(2)} &=-3+\frac{9}{16}t+\frac{3}{64}t^2 f^{-2}(t)-\frac{3}{16}t^{3/2}f^{-1}(t)-\frac{3}{4}\sqrt{t}f(t)+\frac{3}{4}f^2(t),\label{r5:27b}
\end{align}
where
\[ 
f(t)=\left[8+t^{3/2}/8+4\sqrt{4+t^{3/2}/8}\right]^{1/3}.
\]
The next step consists in finding $t=t^{(2)}$ satisfying the condition analogical to   \rf{r5_10}, namely:
\begin{equation}
\forall_{t\neq t^{(2)}}\max_{\tilde{x}\in\tilde{A}}|{\Delta}_{1}^{(2)}(\tilde{x},t^{(2)})|\leq\max_{\tilde{x}\in\tilde{A}}| {\Delta}_{1}^{(2)}(\tilde{x},t)|.     \label{r5:28}
\end{equation}
For this purpose we solve numerically the equation
\begin{equation}
\Delta_1^{(2)}(t,t)=\frac{1}{2}d_1^{(2)}(1+\tilde{\delta}_0(\tilde{x}_0^{II},t))- \frac{1}{2}\tilde{\delta}_{0}^2(\tilde{x}_0^{II},t)(3+\tilde{\delta}_{0}(\tilde{x}_0^{II},t)),\label{r5:29}
\end{equation}
which equates the minimum boundary value of the error of analysed correction with its smallest local minimum, where $x_0^{II}$ is given by \rf{x-max}.
Thus we obtain $t^{(2)}$:
\begin{equation}
t^{(2)}\simeq 3.73157124016,\label{r5-30}
\end{equation}
which is closer to $t_0^{(r)}$ than $t_1^{(r)}$. The value $t=t^{(2)}$ found in this way corresponds to the following magic constant:
\begin{equation}
\tilde{R}^{(2)}=1597466888=0x5F376908            \label{r5-30b}
\end{equation}
and to the relative error 
\begin{equation}
\Delta_{1\max}^{(2)} \simeq 8.7908386407\cdot 10^{-4} \simeq\frac{\tilde{\delta}_{1\max}}{1.99} ,\quad d_1^{(2)}=1.75791023259\cdot 10^{-3} .\label{r5-31}
\end{equation}
In spite of the fact that this error is only slightly greater than $\Delta^{(1)}_{1\max}$, 
\[
\frac{\Delta^{(2)}_{1\max}-\Delta^{(1)}_{1\max}}{\Delta^{(1)}_{1\max}}\simeq 0.3\%,
\]
the difference between error functions, i.e.,
\[
{\Delta_1^{(1)}(\tilde{x},t^{(2)})-\Delta_1^{(2)}(\tilde{x},t^{(1)})} ,
\]
can reach much higher values, even $2\%$ of $\Delta_{1\max}^{(1)}$  (see Fig.~\ref{pic6}), due to a different estimation of the parameter $t$. 

In the case of the second correction, we keep the obtained value $t=t^{(2)}$ and determine the parameter $d_2^{(2)}$  equating the maximum value of the error with the modulus of its global minimum.   ${\Delta^{(2)}_2}(\tilde{x},t^{(2)})$ is increasing (decreasing) with respect to negative (positive) ${\Delta^{(2)}_1}(\tilde{x},t^{(2)})$ and has local minima which come only from positive maxima and negative minima. Therefore the global minimum should correspond to the global minimum $-\Delta_{1\max}^{(2)}$ or to the global maximum $\Delta_{1\max}^{(2)}$. 
Substituting these values to Eq.~\rf{r5:21b} in the place of $\Delta_{1}^{(2)}(\tilde{x},t^{(2)})$ we obtain that deeper minima of $(-\Delta_{2\max}^{(2)})$ come from the global minimum of the first correction:
\begin{equation}
-\Delta_{2\max}^{(2)}=\frac{1}{2}d_2^{(2)}(1 - \Delta_{1\max}^{(2)})- \frac{1}{2}\Delta_{1\max}^{(2)2}(3 - \Delta_{1\max}^{(2)}),
\end{equation}
and the maximum, by analogy to the first correction, corresponds to the following value of   $\Delta_{1}^{(2)}(\tilde{x},t^{(2)})$:  
\begin{equation}
\Delta^{+}=\sqrt{1+d_2^{(2)}/3}-1.
\end{equation}
Solving the equation 
\begin{equation}
\Delta_{2\max}^{(2)}=\frac{1}{2}d_2^{(2)}(1+\Delta^{+})- \frac{1}{2}\Delta^{+2}(3+\Delta^{+})
\end{equation}
we get 
\begin{equation}
d_2^{(2)}\simeq 1.159352515 \cdot 10^{-6}\quad\text{and}\quad \Delta_{2\max}^{(2)}  \simeq =5.796763137\cdot 10^{-7}  \simeq \frac{\tilde{\delta}_{2\max}}{7.93} .
\end{equation} 
Thus we completed the derivation of the function $InvSqrt2$. The computer code contains a new magic constant, see \rf{r5-30b}, and has two lines ($6$ and $7$) modified as compared with the code $InvSqrt$:
\begin{tabbing}
0000000\=122\=5678\=\kill
\>\textit{1.}\> \textbf{float} InvSqrt2(\textbf{float} x)\{ \\
\>\textit{2.}\>\>\textbf{float} halfnumber = 0.5f * x; \\
\>\textit{3.}\>\> \textbf{int} i = *(\textbf{int}*) \&x;\\
\>\textit{4.}\>\> i = 0x5F376908 - (i$>>$1);\\
\>\textit{5.}\>\>  x = *(\textbf{float}*)\&i;\\
\>\textit{6.}\>\> x = x*(1.5008789f - halfnumber*x*x); \\
\>\textit{7.}\>\> x = x*(1.5000006f - halfnumber*x*x); \\
\>\textit{8.}\>\> \textbf{return} x ;\\
\>\textit{9.}\>\}
\end{tabbing}
where
\[
1.5008789f\simeq \frac{3+d_1^{(2)}}{2},\quad\quad 1.5000006f\simeq \frac{3+d_2^{(2)}}{2}.
\]
We point out that the code $InvSqrt2$ has the same number of algebraic operations as $InvSqrt$.

\section{Numerical experiments}

\begin{figure}
\begin{center}
\includegraphics[width=15cm]{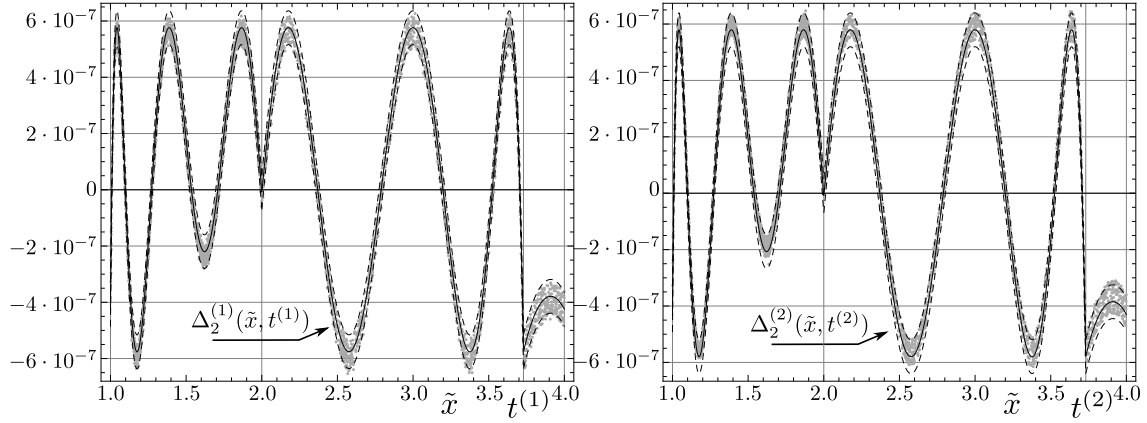}
\end{center}
\caption{Solid lines represent functions $\Delta_2^{(1)}(\tilde{x},t^{(1)})$ (left) and $\Delta_2^{(2)}(\tilde{x},t^{(2)})$ (right). Their vertical shifts by $\pm 6\cdot 10^{-8}$ are denoted by dashed lines. Finally, dots represent relative errors for $4000$ random values $x\in (2^{-126},2^{128})$, produced by algorithms $InvSqrt1$ (left) and $InvSqrt2$ (right).}
\label{pic7}
\end{figure}

\begin{figure}
\begin{center}
\includegraphics[width=14cm]{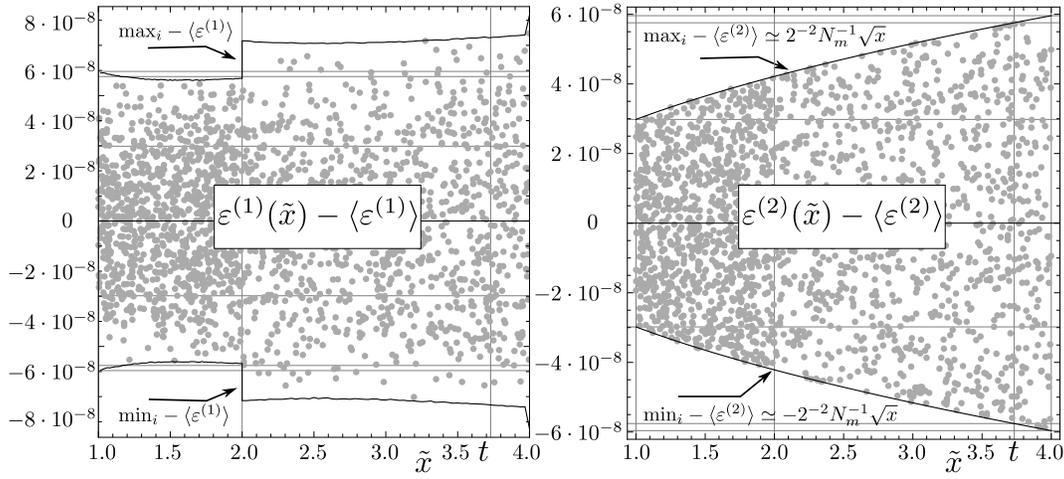}
\end{center}
\caption{Relative errors $\varepsilon^{(1)}$ (left) and  $\varepsilon^{(2)}$ (right) arising during the \textbf{float} approximation of corrections $\tilde{y}_{22}(\tilde{x},t)$ and $\tilde{y}_{32}(\tilde{x},t)$, respectively. Dots represent errors determined for $2000$ random values $\tilde{x}\in [1,4)$. Solid lines represent maximum ($\max_i$) and minimum ($\min_i$) values of relative errors (intervals $[1,2)$ and $[2,4)$ were divided into $64$ equal intervals, and then   extremum values were determined in all these intervals). }
\label{pic8}
\end{figure}

The new algorithms were tested on the processor Intel Core 2 Quad (x86-64) using the compiler  TDM-GCC 4.9.2 32-bit (then, in the case of $InvSqrt$, the values of erors are practically the same as those obtained by Lomont \cite{Lomont}). The same reults were obtained also on Intel i7-5700.

 Applying the algorithm $InvSqrt1$ we obtain relative errors characterized by ``oscillations'' with a center slightly  shifted with respect to the analytical solution, see Fig.~\ref{pic7}. Although at figures we present only randomly chosen values but calculations were carried out for all numbers $x$ of the type \textbf{float} such that $e_x\in[-126,128)$, for any interval $A_n$ (see Eq. (2.7) in \cite{MWHHC}). The range of errors is the same for all these intervals (except $e_x=-126$): 
\begin{equation}  \label{Num1}
\Delta^{(1)}_N(x)=\mathrm{sqrt}(x)*\mathrm{InvSqrt1}(x)-1. \in(\Delta^{(1)}_{N,\min}, \Delta^{(1)}_{N,\max}),
\end{equation}
where 
\[
\Delta^{(1)}_{N,\min}=-6.62\cdot 10^{-7},\quad\quad  \Delta^{(1)}_{N,\max}=6.35\cdot 10^{-7}. 
\]
For $e_x=-126$ the interval of errors is slightly wider: 
\[ [ -6.72\cdot 10^{-7}, 6.49\cdot 10^{-7} ]   .   \] 
The observed blur can be noticed already for the approximation error of the correction $\tilde{y}_{22}(\tilde{x})$:
\begin{equation}
\varepsilon^{(1)}(\tilde{x})= \frac{\mathrm{InvSqrt1(x)}-\tilde{y}_{22}(\tilde{x},t^{(1)})}{\tilde{y}_{22}(\tilde{x},t^{(1)})} .
\end{equation}
The values of this error are distributed symmetrically around the mean value  $\langle\varepsilon^{(1)} \rangle$:
\begin{equation}
\langle\varepsilon^{(1)} \rangle=2^{-1}N_m^{-23}\sum_{x\in [1,4)}\varepsilon^{(1)}(\tilde{x})=-1.398\cdot 10^{-8}
\end{equation}
enclosing the range:
\begin{equation}
\varepsilon^{(1)}(\tilde{x})\in [-9.676\cdot 10^{-8}, 6.805\cdot 10^{-8}] ,
\end{equation}
see Fig.~\ref{pic8}. The blur parameters of the function 
 $\varepsilon^{(1)}(\tilde{x},t)$ show that the main source of the difference between analytical and numerical results is the use of precision \textbf{float} and, in particular, rounding of constant parameters of the function   $\mathrm{InvSqrt1(\tilde{x})}$.  It is worthwhile to point out that in this case the amplitude of the error oscillations is about $40\%$ greater than the amplitude of oscillations of $(\tilde{y}_{00}-\tilde{y}_0)/\tilde{y}_0$ (i.e., in the case of $InvSqrt$), see the right part of  Fig.~2 in \cite{MWHHC}.

The errors of numerical values returned by $InvSqrt2$ belong (for $e_x \neq -126$) to the following interval
\begin{equation}  \label{Num2}
\Delta^{(2)}_N(x)=\mathrm{sqrt}(x)*\mathrm{InvSqrt2}(x)-1.\in(\Delta^{(2)}_{N,\min}, \Delta^{(2)}_{N,\max}),
\end{equation}
where
\[
\Delta^{(2)}_{N,\min}=-6.21\cdot 10^{-7},\quad\quad  \Delta^{(2)}_{N,\max}=6.53\cdot 10^{-7} .
\]
For $e_x = - 126$ we also get a wider interval:
\[ [ -6.46\cdot 10^{-7}, 6.84\cdot 10^{-7} ]   .   \] 
These errors differ  from the errors of  $\tilde{y}_{32}(\tilde{x},t^{(2)})$ determined analytically (like in the case of the function $InvSqrt1$). 
The observed blur of the \textbf{float} approximation of the error of the correction $\tilde{y}_{32}(\tilde{x})$:
\begin{equation}
\varepsilon^{(2)}(\tilde{x})= \frac{\mathrm{InvSqrt2(x)}-\tilde{y}_{32}(\tilde{x},t^{(2)})}{\tilde{y}_{32}(\tilde{x},t^{(2)})},
\end{equation}
is also symmetric with respect to the mean value $\langle\varepsilon^{(2)} \rangle$ (see Fig.~\ref{pic8}):
\begin{equation}
\langle\varepsilon^{(2)} \rangle=2^{-1}N_m^{-23}\sum_{x\in [1,4)}\varepsilon^{(2)}(\tilde{x})=1.653\cdot 10^{-8}
\end{equation}
but it covers a much smaller range of values:
\begin{equation}  
\varepsilon^{(2)}(\tilde{x})\in [-4.316\cdot 10^{-8}, 7.612\cdot 10^{-8} ] \;.
\end{equation}
As a consequence, the maximum numerical error of the function $\mathrm{InvSqrt2}$ is smaller than the error of $\mathrm{InvSqrt1{}}$:
\begin{equation}
\max\{|\Delta_{N\min}^{(1)}|,\,|\Delta_{N\max}^{(1)}|\}=-\Delta_{N\min}^{(1)} > \Delta_{N\max}^{(2)}=\max\{|\Delta_{N\min}^{(2)}|,\,|\Delta_{N\max}^{(2)}|\} .
\end{equation}
Results produced by the same hardware with 64-bit compiler have the amplitude of the error oscillations greater by about $10^{-7}$ as compared with the 32-bit case.  

If we stop both algorithms after the first correction, then the numerical errors of $\mathrm{InvSqrt2}$ are contained in the interval
\[              [-8.79\cdot 10^{-4},\,8.79\cdot 10^{-4}]\,,             \] 
which is a little bit wider than the analogical interval for $\mathrm{InvSqrt1}$
\[              [-8.76\cdot 10^{-4},\,8.76\cdot 10^{-4}]\,.           \]
Therefore, the code $\mathrm{InvSqrt1}$  is slightly more accurate than $\mathrm{InvSqrt2}$ as far as one iteration is concerned and both are about 2 times more accurate than $\mathrm{InvSqrt}$. Comparing $\Delta^{(2)}_{N,\max}$ with the analogical numerical bound for $\mathrm{InvSqrt}$ (see Eq. (4.23) in \cite{MWHHC}) we conclude that in the case of two iterations the code $\mathrm{InvSqrt2}$ is about $7$ times more accurate  than $\mathrm{InvSqrt}$. We point out that round-off errors significantly decrease the theoretical improvement which is given by the factor $8$.

\section{Conclusions}

In this paper we have presented two modifications of the famous \textit{InvSqrt} code for fast computation of the inverse square root. The first modification, denoted by \textit{InvSqrt1}, has the same magic constant as \textit{InvSqrt} (i.e., the same initial values for the iteration process) but instead of the standard Newton-Raphson method we propose a similar procedure but  with different coefficients. The obtained algorithm is slightly more expensive in the case of two iterations (8 multiplications instead of 7 at every step) but much more accurate than $InvSqrt$: two times more accurate after the first iteration and about 7  times more accurate after two iterations. The second modification,  denoted by \textit{InvSqrt2}, uses another magic constant. Its  computational cost is practically the same as the cost of the \textit{InvSqrt} code (both have the same number of multiplications). This is the main advantage of \textit{InvSqrt2} over \textit{InvSqrt1} because the accuracy of  \textit{InvSqrt2} is similar to the accuracy of \textit{InvSqrt1}. Analytical values for \textit{InvSqrt1} are slightly better but numerical tests show that round-off errors are a little bit greater in this case,  see \rf{Num1} and \rf{Num2}, although the situation becomes reversed  when the algorithms are stopped after only one iteration. 

Concerning potential applications, we have to acknowlede that for general purpose computing the SSE reciprocal square root instruction is faster and more accurate. We hope, however, that the proposed algorithms  can be applied in embedded systems and microcontrollers that lack a FPU, and potentially in FPGA's.

\end{document}